\documentclass{article}

\usepackage{amsmath,amssymb,amsthm}
\usepackage{graphicx}

\newtheorem{theorem}{Theorem}
\newtheorem{proposition}[theorem]{Proposition}
\newtheorem{lemma}[theorem]{Lemma}
\newtheorem{corollary}[theorem]{Corollary}
\newtheorem{example}{Example}
\newtheorem{definition}[theorem]{Definition}
\newtheorem{remark}{Remark}

\begin{document}

\title{The delta-nabla calculus of variations\thanks{Accepted for publication (02/December/2009) in \emph{Fasciculi Mathematici}.}}

\author{Agnieszka B. Malinowska$^{1, 2}$\\
\texttt{abmalinowska@ua.pt}
\and
Delfim F. M. Torres$^{1}$\\
\texttt{delfim@ua.pt}}

\date{$^1$Department of Mathematics\\
      University of Aveiro\\
      3810-193 Aveiro, Portugal\\[0.3cm]
      $^2$Faculty of Computer Science\\
      Bia{\l}ystok University of Technology\\
      15-351 Bia\l ystok, Poland}

\maketitle


\begin{abstract}
The discrete-time, the quantum,
and the continuous calculus of variations
have been recently unified and extended.
Two approaches are followed in the literature: one dealing with minimization
of delta integrals; the other dealing with minimization
of nabla integrals. Here we propose a more general approach to the
calculus of variations on time scales that allows to obtain
both delta and nabla results as particular cases.

\bigskip


\medskip

\end{abstract}


\section{Introduction}

The calculus of variations on time scales was introduced
by M.~Bohner using the delta derivative and integral
\cite{B:04}: to extremize a functional of the form
\begin{equation}
\label{eq:Pd}
\mathcal{J}_\Delta(y)
= \int_a^b L\left(t,y^\sigma(t),y^\Delta(t)\right) \Delta t \, .
\end{equation}
Motivated by applications in economics \cite{A:B:L:06,A:U:08},
a different formulation for the problems of the calculus
of variations on time scales has been considered, which
involve a functional with a nabla derivative and
a nabla integral \cite{A:T,Atici:comparison,NM:T}:
\begin{equation}
\label{eq:Pn}
\mathcal{J}_\nabla(y)
= \int_a^b L\left(t,y^\rho(t),y^\nabla(t)\right) \nabla t \, .
\end{equation}
Formulations \eqref{eq:Pd} and \eqref{eq:Pn} are consistent
in the sense that results obtained \emph{via}
delta and nabla approaches are similar among them
and similar to the classical results of the calculus
of variations. An example of this is given
by the time scale versions of the Euler-Lagrange equations:
if $y \in C_{\textrm{rd}}^2$ is an extremizer of \eqref{eq:Pd},
then $y$ satisfies the delta-differential equation
\begin{equation}
\label{eq:el:d} \frac{\Delta}{\Delta t}
\partial_{3}L\left(t,y^\sigma(t),{y}^\Delta(t)\right)
= \partial_{2}L\left(t,y^\sigma(t),{y}^\Delta(t)\right)
\end{equation}
for all $t \in [a,b]^{\kappa^2}$ \cite{B:04};
if $y \in C_{\textrm{ld}}^2$ is an extremizer of \eqref{eq:Pn},
then $y$ satisfies the nabla-differential equation
\begin{equation}
\label{eq:el:n} \frac{\nabla}{\nabla t}
\partial_{3}L\left(t,y^\rho(t),{y}^\nabla(t)\right)
= \partial_{2}L\left(t,y^\rho(t),{y}^\nabla(t)\right)
\end{equation}
for all $t \in [a,b]_{\kappa^2}$ \cite{NM:T},
where we use $\partial_{i}L$ to denote the standard
partial derivative of $L(\cdot,\cdot,\cdot)$
with respect to its $i$th variable, $i = 1,2,3$.
In the classical context $\mathbb{T} = \mathbb{R}$
one has
\begin{equation}
\label{eq:Pc}
\mathcal{J}_\Delta(y)
= \mathcal{J}_\nabla(y)
= \int_a^b L\left(t,y(t),y'(t)\right) dt
\end{equation}
and both \eqref{eq:el:d} and \eqref{eq:el:n} coincide
with the standard Euler-Lagrange equation:
if $y \in C^2$ is an extremizer of
the integral functional \eqref{eq:Pc}, then
\begin{equation*}
\frac{d}{d t}
\partial_{3}L\left(t,y(t),y'(t)\right)
=\partial_{2}L\left(t,y(t),y'(t)\right)
\end{equation*}
for all $t \in [a,b]$.
However, the problems of extremizing \eqref{eq:Pd}
and \eqref{eq:Pn} are intrinsically different, in the sense
that is not possible to obtain the nabla results as corollaries
of the delta ones and \emph{vice versa}. Indeed,
if admissible functions $y$ are of class $C^2$ then (\textrm{cf.} \cite{G:G:S:05})
\begin{equation*}
\mathcal{J}_\Delta(y)
= \int_a^b L\left(t,y^\sigma(t),y^\Delta(t)\right) \Delta t
= \int_a^b L\left(\rho(t),(y^\sigma)^\rho(t),y^\nabla(t)\right) \nabla t
\end{equation*}
while
\begin{equation*}
\mathcal{J}_\nabla(y)
= \int_a^b L\left(t,y^\rho(t),y^\nabla(t)\right) \nabla t
= \int_a^b L\left(\sigma(t),(y^\rho)^\sigma(t),y^\Delta(t)\right) \Delta t
\end{equation*}
and one easily see that functionals
\eqref{eq:Pd} and \eqref{eq:Pn}
have a different nature and are not compatible with each other.
In this paper we introduce a more general formulation
of the calculus of variations that
includes, as trivial examples, the problems with functionals
$\mathcal{J}_\Delta(y)$
and $\mathcal{J}_\nabla(y)$ that have been
previously studied in the literature.
Our main result provides an Euler-Lagrange
necessary optimality type condition
(\textrm{cf.} Theorem~\ref{thm:mr}).


\section{Our goal}

Let $\mathbb{T}$ be a given time scale with $a, b \in \mathbb{T}$, $a < b$,
and $\left(\mathbb{T} \setminus \{a,b\}\right)\cap [a,b] \ne \emptyset$;
$L_{\Delta}(\cdot,\cdot,\cdot)$ and $L_{\nabla}(\cdot,\cdot,\cdot)$ be two given smooth
functions from $\mathbb{T} \times \mathbb{R}^2$ to $\mathbb{R}$.
The results here discussed are trivially generalized for
admissible functions $y : \mathbb{T}\rightarrow\mathbb{R}^n$
but for simplicity of presentation
we restrict ourselves to the scalar case $n=1$.
We consider the delta-nabla integral functional
\begin{equation}
\label{eq:P}
\begin{split}
\mathcal{J}(y)
&= \left(\int_a^b L_{\Delta}\left(t,y^\sigma(t),y^\Delta(t)\right) \Delta t\right)
\cdot
\left(\int_a^b L_{\nabla}\left(t,y^\rho(t),y^\nabla(t)\right) \nabla t\right)\\
&= \int_a^b \int_a^b \left[
L_{\Delta}\left(t,y^\sigma(t),y^\Delta(t)\right)
\cdot L_{\nabla}\left(\tau,y^\rho(\tau),y^\nabla(\tau)\right)
\right] \Delta t \nabla \tau \, .
\end{split}
\end{equation}

\begin{remark}
\label{obs}
In the particular case $L_\nabla \equiv \frac{1}{b-a}$
functional \eqref{eq:P} reduces to \eqref{eq:Pd}
(\textrm{i.e.}, $\mathcal{J}(y) = \mathcal{J}_\Delta(y)$);
in the particular case $L_\Delta \equiv \frac{1}{b-a}$
functional \eqref{eq:P} reduces to \eqref{eq:Pn}
(\textrm{i.e.}, $\mathcal{J}(y) = \mathcal{J}_\nabla(y)$).
\end{remark}

\medskip

Our main goal is to answer the following question:
\emph{What is the Euler-Lagrange equation for
$\mathcal{J}(y)$ defined by \eqref{eq:P}?}

\medskip

For simplicity of notation we introduce the operators
$[y]$ and $\{y\}$ defined by
$[y](t) = \left(t,y^\sigma(t),y^\Delta(t)\right)$ and
$\{y\}(t) = \left(t,y^\rho(t),y^\nabla(t)\right)$.
Then,
\begin{gather*}
\mathcal{J}_\Delta(y) = \int_a^b L_\Delta[y](t) \Delta t \, , \quad
\mathcal{J}_\nabla(y) = \int_a^b L_\nabla\{y\}(t) \nabla t \, , \\
\mathcal{J}(y) = \mathcal{J}_\Delta(y) \mathcal{J}_\nabla(y)
= \int_a^b \int_a^b
L_{\Delta}[y](t) L_{\nabla}\{y\}(\tau) \Delta t \nabla \tau \, .
\end{gather*}


\section{Preliminaries to the calculus of variations}

Similar to the classical calculus of variations,
integration by parts will play an important role in our delta-nabla
calculus of variations.
If functions $f,g : \mathbb{T}\rightarrow\mathbb{R}$
are delta and nabla differentiable with continuous derivatives,
then the following formulas of integration by parts hold \cite{B:P:01}:
\begin{equation}
\label{intBP}
\begin{split}
\int_{a}^{b}f^\sigma(t) g^{\Delta}(t)\Delta t
&=\left.(fg)(t)\right|_{t=a}^{t=b}
-\int_{a}^{b}f^{\Delta}(t)g(t)\Delta t \, , \\
\int_{a}^{b}f(t)g^{\Delta}(t)\Delta t
&=\left.(fg)(t)\right|_{t=a}^{t=b}
-\int_{a}^{b}f^{\Delta}(t)g^\sigma(t)\Delta t \, , \\
\int_{a}^{b}f^\rho(t)g^{\nabla}(t)\nabla t
&=\left.(fg)(t)\right|_{t=a}^{t=b}
-\int_{a}^{b}f^{\nabla}(t)g(t)\nabla t \, ,\\
\int_{a}^{b}f(t)g^{\nabla}(t)\nabla t
&=\left.(fg)(t)\right|_{t=a}^{t=b}
-\int_{a}^{b}f^{\nabla}(t)g^\rho(t)\nabla t \, .
\end{split}
\end{equation}

The following fundamental lemma of the calculus of variations
on time scales involving a nabla derivative
and a nabla integral has been proved in \cite{NM:T}.

\begin{lemma}{\rm (The nabla Dubois-Reymond lemma \cite[Lemma~14]{NM:T}).}
\label{DBRL:n}
Let $f \in C_{\textrm{ld}}([a,b], \mathbb{R})$. If
$$
\int_{a}^{b} f(t)\eta^{\nabla}(t)\nabla t=0 \quad
\mbox{for all $\eta \in C_{\textrm{ld}}^1([a,b],
\mathbb{R})$ with $\eta(a)=\eta(b)=0$} \, ,
$$
then $f(t) = c$ on $t\in [a,b]_\kappa$
for some constant $c$.
\end{lemma}

Lemma~\ref{DBRL:d} is the analogous delta version of Lemma~\ref{DBRL:n}:

\begin{lemma}{\rm (The delta Dubois-Reymond lemma \cite{B:04}).}
\label{DBRL:d}
Let $g\in C_{\textrm{rd}}([a,b], \mathbb{R})$. If
$$\int_{a}^{b} g(t) \eta^\Delta(t)\Delta t=0  \quad
\mbox{for all $\eta \in C_{\textrm{rd}}^1$ with
$\eta(a)=\eta(b)=0$} \, ,$$
then $g(t)=c$ on $[a,b]^\kappa$
for some $c\in\mathbb{R}$.
\end{lemma}

Proposition~\ref{prop:rel:der}
gives a relationship between delta
and nabla derivatives.

\begin{proposition}{\rm (Theorems~2.5 and 2.6 of \cite{A:G:02}).}
\label{prop:rel:der}
(i) If $f : \mathbb{T} \rightarrow \mathbb{R}$ is delta differentiable
on $\mathbb{T}^\kappa$ and $f^\Delta$ is continuous on $\mathbb{T}^\kappa$,
then $f$ is nabla differentiable on $\mathbb{T}_\kappa$ and
\begin{equation}
\label{eq:chgN_to_D}
f^\nabla(t)=\left(f^\Delta\right)^\rho(t) \quad \text{for all }
t \in \mathbb{T}_\kappa \, .
\end{equation}
(ii) If $f : \mathbb{T} \rightarrow \mathbb{R}$ is nabla differentiable
on $\mathbb{T}_\kappa$ and $f^\nabla$ is continuous on $\mathbb{T}_\kappa$,
then $f$ is delta differentiable on $\mathbb{T}^\kappa$ and
\begin{equation}
\label{eq:chgD_to_N}
f^\Delta(t)=\left(f^\nabla\right)^\sigma(t) \quad \text{for all }
t \in \mathbb{T}^\kappa \, .
\end{equation}
\end{proposition}

\begin{remark}
Note that, in general, $f^\nabla(t) \ne f^\Delta\left(\rho(t)\right)$
and $f^\Delta(t) \ne f^\nabla\left(\sigma(t)\right)$.
In Proposition~\ref{prop:rel:der} the assumptions on the
continuity of $f^\Delta$ and $f^\nabla$ are crucial.
\end{remark}

\begin{proposition}{\rm (\cite[Theorem~2.8]{A:G:02}).}
\label{eq:prop}
Let $a, b \in\mathbb{T}$ with $a \le b$ and let $f$
be a continuous function on $[a, b]$. Then,
\begin{equation*}
\begin{split}
\int_a^b f(t)\Delta t
&= \int_a^{\rho(b)} f(t)\Delta t
+ (b - \rho(b))f^\rho(b) \, , \\
\int_a^b f(t)\Delta t
&= (\sigma(a) - a) f(a)
+ \int_{\sigma(a)}^b f(t)\Delta t \, , \\
\int_a^b f(t)\nabla t
&= \int_a^{\rho(b)} f(t)\nabla t
+ (b - \rho(b)) f(b) \, , \\
\int_a^b f(t)\nabla t
&= (\sigma(a) - a) f^\sigma(a)
+ \int_{\sigma(a)}^b f(t)\nabla t \, .
\end{split}
\end{equation*}
\end{proposition}

We end our brief review of the calculus on time scales
with a relationship between the delta and nabla integrals.

\begin{proposition}{\rm (\cite[Proposition~7]{G:G:S:05}).}
If function $f : \mathbb{T} \rightarrow \mathbb{R}$
is continuous, then for all $a, b \in \mathbb{T}$
with $a < b$ we have
\begin{gather}
\int_a^b f(t) \Delta t = \int_a^b f^\rho(t) \nabla t \, , \label{eq:DtoN}\\
\int_a^b f(t) \nabla t = \int_a^b f^\sigma(t) \Delta t \, . \label{eq:NtoD}
\end{gather}
\end{proposition}


\section{Main Result}

We consider the problem of extremizing
the variational functional \eqref{eq:P}
subject to given boundary conditions
$y(a) = \alpha$ and $y(b) = \beta$:
\begin{equation}
\label{problem:P}
\begin{gathered}
\mathcal{J}(y) =
\left(\int_a^b L_{\Delta}[y](t) \Delta t\right)
\left(\int_a^b L_{\nabla}\{y\}(t) \nabla t\right) \longrightarrow
\textrm{extr} \\
y(\cdot) \in C_{\diamond}^1 \\
y(a) = \alpha \, , \quad y(b) = \beta \, ,
\end{gathered}
\end{equation}
where $C_{\diamond}^1$ denote the class of functions
$y : [a,b]\rightarrow\mathbb{R}$  with
$y^\Delta$ continuous on $[a,b]^\kappa$
and $y^\nabla$ continuous on $[a,b]_\kappa$.
Before presenting the Euler-Lagrange equations for
problem \eqref{problem:P} we introduce the definition
of weak local extremum.

\begin{definition} We say that $\hat{y}\in C_{\diamond}^{1}([a,b],
\mathbb{R})$ is a weak local minimizer (respectively weak local
maximizer) for problem \eqref{problem:P} if there exists
$\delta >0$ such that
$\mathcal{J}(\hat{y})\leq \mathcal{J}(y)$
(respectively $\mathcal{J}(\hat{y})
\geq \mathcal{J}(y)$) for all
$y \in C_{\diamond}^{1}([a,b], \mathbb{R})$ satisfying the boundary
conditions $y(a) = \alpha$, $y(b) = \beta$, and
$\parallel y - \hat{y}\parallel_{1,\infty} < \delta$,
where
$$\parallel y\parallel_{1,\infty}:=
\parallel y^{\sigma}\parallel_{\infty}
+ \parallel y^{\rho}\parallel_{\infty}
+ \parallel y^{\Delta}\parallel_{\infty}
+ \parallel y^{\nabla}\parallel_{\infty}$$
and $\parallel y\parallel_{\infty} :=\sup_{t \in
[a,b]_{\kappa}^{\kappa}}\mid y(t) \mid$.
\end{definition}


Theorem~\ref{thm:mr} gives two different forms
for the Euler-Lagrange equation on time scales
associated with the variational problem \eqref{problem:P}.

\begin{theorem}{\rm (The general Euler-Lagrange equations
on time scales).}
\label{thm:mr}
If $\hat{y} \in C_{\diamond}^1$ is a weak local extremizer of problem
\eqref{problem:P}, then $\hat{y}$ satisfies
the following delta-nabla integral equations:
\begin{multline}
\label{eq:EL1}
\mathcal{J}_\nabla(\hat{y})
\left(\partial_3 L_\Delta[\hat{y}](\rho(t))
-\int_{a}^{\rho(t)} \partial_2 L_\Delta[\hat{y}](\tau) \Delta\tau\right)\\
+
\mathcal{J}_\Delta(\hat{y})
\left(\partial_3 L_\nabla\{\hat{y}\}(t)
-\int_{a}^{t} \partial_2 L_\nabla\{\hat{y}\}(\tau) \nabla\tau\right)
= \text{const} \quad \forall t \in [a,b]_\kappa \, ;
\end{multline}
\begin{multline}
\label{eq:EL2}
\mathcal{J}_\nabla(\hat{y})
\left(\partial_3 L_\Delta[\hat{y}](t)
-\int_{a}^{t} \partial_2 L_\Delta[\hat{y}](\tau) \Delta\tau\right)\\
+
\mathcal{J}_\Delta(\hat{y})
\left(\partial_3 L_\nabla\{\hat{y}\}(\sigma(t))
-\int_{a}^{\sigma(t)} \partial_2 L_\nabla\{\hat{y}\}(\tau) \nabla\tau\right)
= \text{const} \quad \forall t \in [a,b]^\kappa \, .
\end{multline}
\end{theorem}

\begin{remark}
In the classical context (\textrm{i.e.},
when $\mathbb{T} = \mathbb{R}$) the
necessary conditions \eqref{eq:EL1}
and \eqref{eq:EL2} coincide with the Euler-Lagrange
equations recently given in \cite{Pedregal}.
\end{remark}

\begin{proof}
Suppose that $\mathcal{J}$
has a weak local extremum at $\hat{y}$. We
consider the value of $\mathcal{J}$ at nearby
functions $\hat{y} + \varepsilon \eta$,
where $\varepsilon\in \mathbb{R}$ is a small parameter,
$\eta \in C_{\diamond}^{1}([a,b],\mathbb{R})$ with $\eta(a)=\eta(b)=0$.
Thus, function $\phi(\varepsilon)
= \mathcal{J}(\hat{y} + \varepsilon \eta)$
has an extremum at $\varepsilon = 0$. Using the first-order necessary
optimality condition $\left.\phi'(\varepsilon)\right|_{\varepsilon = 0} = 0$,
\begin{multline}
\label{eq:prf:+}
\mathcal{J}_\nabla(\hat{y}) \int_a^b
\left(\partial_2 L_\Delta[\hat{y}](t) \eta^\sigma(t)
+ \partial_3 L_\Delta[\hat{y}](t) \eta^\Delta(t)\right) \Delta t
\\
+ \mathcal{J}_\Delta(\hat{y}) \int_a^b
\left(\partial_2 L_\nabla\{\hat{y}\}(t) \eta^\rho(t)
+ \partial_3 L_\nabla\{\hat{y}\}(t) \eta^\nabla(t)\right)
\nabla t = 0 \, .
\end{multline}
Let $A(t) = \int_a^t \partial_2 L_\Delta[\hat{y}](\tau) \Delta\tau$
and $B(t) = \int_a^t \partial_2 L_\nabla\{\hat{y}\}(\tau) \nabla\tau$.
Then, $A^\Delta(t) = \partial_2 L_\Delta[\hat{y}](t)$,
$B^\nabla(t) = \partial_2 L_\nabla\{\hat{y}\}(t)$,
and the first and third integration by parts formula
in \eqref{intBP} tell us, respectively, that
\begin{equation*}
\begin{split}
\int_a^b \partial_2 L_\Delta[\hat{y}](t) \eta^\sigma(t) \Delta t
&= \int_a^b A^\Delta(t) \eta^\sigma(t) \Delta t
= \left. A(t) \eta(t)\right|_{t=a}^{t=b} - \int_a^b A(t) \eta^\Delta(t) \Delta t\\
&= - \int_a^b A(t) \eta^\Delta(t) \Delta t
\end{split}
\end{equation*}
and
\begin{equation*}
\begin{split}
\int_a^b \partial_2 L_\nabla\{\hat{y}\}(t) \eta^\rho(t) \nabla t
&= \int_a^b B^\nabla(t) \eta^\rho(t) \nabla t
= \left. B(t) \eta(t)\right|_{t=a}^{t=b}
- \int_a^b B(t) \eta^\nabla(t) \nabla t\\
&= - \int_a^b B(t) \eta^\nabla(t) \nabla t \, .
\end{split}
\end{equation*}
If we denote $f(t) = \partial_3 L_\Delta[\hat{y}](t) - A(t)$
and $g(t) = \partial_3 L_\nabla\{\hat{y}\}(t) - B(t)$,
then we can write the necessary optimality condition
\eqref{eq:prf:+} in the form
\begin{equation}
\label{eq:prf:+:aftIP}
\mathcal{J}_\nabla(\hat{y}) \int_a^b f(t) \eta^\Delta(t) \Delta t
+ \mathcal{J}_\Delta(\hat{y}) \int_a^b g(t) \eta^\nabla(t) \nabla t = 0 \, .
\end{equation}
We now split the proof in two parts:
we prove \eqref{eq:EL1} transforming the delta integral
in \eqref{eq:prf:+:aftIP} to a nabla integral by means of
\eqref{eq:DtoN}; we prove \eqref{eq:EL2} transforming
the nabla integral in \eqref{eq:prf:+:aftIP}
to a delta integral by means of \eqref{eq:NtoD}.
By \eqref{eq:DtoN} the necessary optimality
condition \eqref{eq:prf:+:aftIP} is equivalent to
\begin{equation*}
\int_a^b \left(\mathcal{J}_\nabla(\hat{y})
f^\rho(t) (\eta^\Delta)^\rho(t) + \mathcal{J}_\Delta(\hat{y})
g(t) \eta^\nabla(t)\right) \nabla t = 0
\end{equation*}
and by \eqref{eq:chgN_to_D} to
\begin{equation}
\label{eq:bef:FL1}
\int_a^b \left(\mathcal{J}_\nabla(\hat{y}) f^\rho(t)
+ \mathcal{J}_\Delta(\hat{y}) g(t)\right)
\eta^\nabla(t) \nabla t = 0 \, .
\end{equation}
Applying Lemma~\ref{DBRL:n} to \eqref{eq:bef:FL1}
we prove \eqref{eq:EL1}:
\begin{equation*}
\mathcal{J}_\nabla(\hat{y}) f^\rho(t)
+ \mathcal{J}_\Delta(\hat{y}) g(t) = c \quad \forall t \in [a,b]_\kappa \, ,
\end{equation*}
where $c$ is a constant. By \eqref{eq:NtoD} the necessary optimality
condition \eqref{eq:prf:+:aftIP} is equivalent to
$\int_a^b \left(\mathcal{J}_\nabla(\hat{y}) f(t) \eta^\Delta(t)
+ \mathcal{J}_\Delta(\hat{y}) g^\sigma(t)
\left(\eta^\nabla\right)^\sigma(t)\right) \Delta t = 0$
and by \eqref{eq:chgD_to_N} to
\begin{equation}
\label{eq:bef:FL2}
\int_a^b \left(\mathcal{J}_\nabla(\hat{y}) f(t)
+ \mathcal{J}_\Delta(\hat{y}) g^\sigma(t)\right)
\eta^\Delta(t) \Delta t = 0 \, .
\end{equation}
Applying Lemma~\ref{DBRL:d} to \eqref{eq:bef:FL2}
we prove \eqref{eq:EL2}:
\begin{equation*}
\mathcal{J}_\nabla(\hat{y}) f(t)
+ \mathcal{J}_\Delta(\hat{y}) g^\sigma(t)
= c \quad \forall t \in [a,b]^\kappa \, ,
\end{equation*}
where $c$ is a constant.
\end{proof}

\begin{corollary}
Let $L_\Delta\left(t,y^\sigma,y^\Delta\right) = L_\Delta(t)$
and $\mathcal{J}_\Delta(\hat{y}) \ne 0$
(this is true, \textrm{e.g.}, for
$L_\Delta \equiv \frac{1}{b-a}$ for which $\mathcal{J}_\Delta
= 1$; \textrm{cf.} Remark~\ref{obs}). Then, $\partial_2 L_\Delta
= \partial_3 L_\Delta = 0$ and the Euler-Lagrange equation \eqref{eq:EL1} takes the form
\begin{equation}
\label{cor:EL1}
\partial_3 L_\nabla\{\hat{y}\}(t)
-\int_{a}^{t} \partial_2 L_\nabla\{\hat{y}\}(\tau) \nabla\tau
= \text{const} \quad \forall t \in [a,b]_\kappa \, .
\end{equation}
\end{corollary}

\begin{remark}
If $\hat{y} \in C_{\textrm{ld}}^2$, then
nabla-differentiating \eqref{cor:EL1} we obtain
the Euler-Lagrange differential equation
\eqref{eq:el:n} as proved in \cite{NM:T}:
\begin{equation*}
\frac{\nabla}{\nabla t}
\partial_3 L_\nabla\{\hat{y}\}(t)
- \partial_2 L_\nabla\{\hat{y}\}(t) = 0 \quad \forall t \in [a,b]_{\kappa^2} \, .
\end{equation*}
\end{remark}

\begin{corollary}
Let $L_\nabla\left(t,y^\rho,y^\nabla\right) = L_\nabla(t)$
and $\mathcal{J}_\nabla(\hat{y}) \ne 0$
(this is true, \textrm{e.g.}, for
$L_\nabla \equiv \frac{1}{b-a}$ for which $\mathcal{J}_\nabla
= 1$; \textrm{cf.} Remark~\ref{obs}). Then, $\partial_2 L_\nabla
= \partial_3 L_\nabla = 0$ and the Euler-Lagrange equation \eqref{eq:EL2} takes the form
\begin{equation}
\label{cor:EL2}
\partial_3 L_\Delta[\hat{y}](t)
-\int_{a}^{t} \partial_2 L_\Delta[\hat{y}](\tau) \Delta\tau
= \text{const} \quad \forall t \in [a,b]^\kappa \, .
\end{equation}
\end{corollary}

\begin{remark}
If $\hat{y} \in C_{\textrm{rd}}^2$,
then delta-differentiating \eqref{cor:EL2} we obtain the Euler-Lagrange
differential equation \eqref{eq:el:d} as proved in \cite{B:04}:
\begin{equation*}
\frac{\Delta}{\Delta t}
\partial_3 L_\Delta[\hat{y}](t)
- \partial_2 L_\Delta[\hat{y}](t) = 0 \quad \forall t \in [a,b]^{\kappa^2} \, .
\end{equation*}
\end{remark}

\begin{example}
\label{ex:first:simp:ex}
Let $\mathbb{T}$ be a time scale with $0$, $\xi \in \mathbb{T}$,
$0 < \xi$, and $\left(\mathbb{T}\setminus \{0,\xi\}\right)\cap [0,\xi] \ne \emptyset$.
Consider the problem
\begin{equation}
\label{ex:1}
\begin{gathered}
\textrm{minimize} \quad
\mathcal{J}(y)=\left(\int_{0}^{\xi}(y^\Delta(t))^2\Delta t\right)
\left(\int_{0}^{\xi}\left(y^\nabla(t))^2\right)\nabla t\right) \, ,\\
y(0)=0, \quad y(\xi)=\xi \, .
\end{gathered}
\end{equation}
Since $L_{\Delta}=(y^\Delta)^2$ and $L_{\nabla}=(y^\nabla)^2$,
we have $\partial_2L_{\Delta}=0$, $\partial_3L_{\Delta}=2y^\Delta$,
$\partial_2L_{\nabla}=0$, and $\partial_3L_{\nabla}=2y^\nabla$.
Using equation \eqref{eq:EL2} of Theorem~\ref{thm:mr}
we get the following delta-nabla differential equation:
\begin{equation}\label{ex:2}
2Ay^{\Delta}(t)+2By^{\nabla}(\sigma(t))=C,
\end{equation}
where $C\in\mathbb{R}$ and $A$, $B$ are the values of functionals
$\mathcal{J}_{\nabla}$ and $\mathcal{J}_{\Delta}$ in a solution of
problem \eqref{ex:1}, respectively. From \eqref{eq:chgD_to_N}
we can write equation \eqref{ex:2} in the form
\begin{equation}
\label{ex:3}
2Ay^{\Delta}(t)+2By^\Delta=C.
\end{equation}
Observe that $A+B$ cannot be equal to $0$. Thus, solving equation
\eqref{ex:3} subject to the boundary conditions $y(0)=0$ and
$y(\xi)=\xi$ we get $y(t)=t$ as a candidate local minimizer for the
problem \eqref{ex:1}.
\end{example}


\section{Conclusion}

A general necessary optimality condition
for problems of the calculus of variations
on time scales has been given.
The proposed calculus of variations
extends the problems with delta derivatives
considered in \cite{B:T:08,B:04} and analogous nabla
problems \cite{A:T,NM:T} to more general cases
described by the product of a delta and a nabla
integral. Minimization of functionals
given by the product of two integrals
were considered by Euler himself,
and are now receiving an increasing interest
because of their nonlocal properties
and their applications in economics \cite{Pedregal}.


\section*{Acknowledgments}

Work supported by the {\it Centre for Research on Optimization
and Control} (CEOC) from the ``Funda\c{c}\~{a}o
para a Ci\^{e}ncia e a Tecnologia'' (FCT),
cofinanced by the European Community Fund
FEDER/POCI 2010. Agnieszka Malinowska is also supported by
Bia{\l}ystok University of Technology,
via a project of the Polish
Ministry of Science and Higher Education ``Wsparcie miedzynarodowej mobilnosci naukowcow''.



\end{document}